\documentclass[fleqn]{mat01}
\usepackage{times,mathtimy,amssymb,latexsym,amscd}
\begin{document}

\setcounter{page}{399}
\firstpage{399}

\newtheorem{pot}[defin]{Proof of Theorem}{\it }{\rm }
\renewcommand\theequation{\arabic{section}.\arabic{equation}}

\def\rr{{\Bbb R}}
\def\rn{{{\rr}^n}}

\def\dsum{\displaystyle\sum}
\def\dint{\displaystyle\int}
\def\dfrac{\displaystyle\frac}
\def\dsup{\displaystyle\sup}
\def\sumj{\sum_{j=-\infty}^{\infty}}
\def\dsumj{\dsum_{j=-\infty}^{\infty}}

\def\T{T_{\Omega, \alpha}}
\def\tT{T_{\tOmega, \alpha}}
\def\L{{\rm Lip}_{\beta}(\rn)}
\def\tOmega{\tilde \Omega}
\def\tomega{\tilde \omega}

\newcommand{\A}{\mathbb{A}}
\newcommand{\R}{\mathbb{R}}
\newcommand{\reals}{\mathbb{R}}
\newcommand{\C}{\mathbb{C}}
\newcommand{\N}{\mathbb{N}}
\newcommand{\D}{\mathbb{D}}
\newcommand{\Z}{\mathbb{Z}}
\newcommand{\integers}{\mathbb{Z}}
\newcommand{\Q}{\mathbb{Q}}
\newcommand{\HH}{\mathbb{H}}
\newcommand{\K}{\mathbb{K}}
\newcommand{\E}{\mathbb{E}}
\newcommand{\F}{\mathbb{F}}

\newtheorem{defi}[defin]{\rm DEFINITION}

\newtheorem{theor}{\bf Theorem}
\renewcommand\thetheor{\Alph{theor}}

\title{Commutators of integral operators with variable kernels\\ on
Hardy spaces}

\markboth{Pu Zhang and Kai Zhao}{Commutators of integral operators with variable kernels}

\author{PU ZHANG$^{1,*}$ and KAI ZHAO$^{2}$}

\address{$^{1}$Institute of Mathematics, Zhejiang Sci-Tech University,
Hangzhou 310018, People's Republic of China\\
\noindent$^{2}$Department of Mathematics, Qingdao University,
Qingdao 266071, People's Republic of China\\
\noindent $^{*}$Corresponding author\\
\noindent E-mail: puzhang@sohu.com; zhkai01@sina.com}

\volume{115}

\mon{November}

\parts{4}

\pubyear{2005}

\Date{MS received 26 January 2005; revised 13 July 2005}\vspace{-.6pc}

\begin{abstract}
Let $\T \,(0\leq \alpha <n)$ be the singular and fractional
integrals with variable kernel $\Omega(x,z)$, and $[b,\T]$ be the
commutator generated by $\T$ and a Lipschitz function $b$. In this
paper, the authors study the boundedness of $[b,\T]$ on the Hardy
spaces, under some assumptions such as the $L^r$-Dini condition.
Similar results and the weak type estimates at the end-point cases
are also given for the homogeneous convolution operators $\tT
\,(0\leq \alpha <n)$. The smoothness conditions imposed on
$\tOmega$ are weaker than the corresponding known results.\vspace{-.3pc}
\end{abstract}

\keyword{Singular and fractional integrals; variable kernel;
commutator; Hardy space.\vspace{-.3pc}}

\maketitle

\section{Introduction and main results}

Let $S^{n-1}$ be the unit sphere in $\rn(n \geq 2)$ equipped with
the normalized Lebesgue measure $\hbox{d}\sigma (x')$. We say a
function $\Omega(x,z)$ defined on $\rn \times \rn$ belongs to
$L^\infty(\rn)\times L^r(S^{n-1})\ (r\geq 1)$, if
\begin{enumerate}
\renewcommand{\labelenumi}{(\roman{enumi})}
\leftskip .1pc
\item For all $x, z\in \rn$ and $\lambda>0$,
$\Omega(x,\lambda z)=\Omega(x,z);$
\item $\|\Omega\|_{L^\infty(\rn)\times L^r(S^{n-1})}:=
\sup_{x\in \rn}\left(\int_{S^{n-1}}|\Omega(x,z')|^r
\hbox{d}\sigma(z')\right)^{1/r}<\infty.$
\end{enumerate}

For $0\leq \alpha<n,$ we define the integral operators with
variable kernels as follows:
\begin{equation}
T_{\Omega,\alpha}f(x)=\int_{\rn}\frac {\Omega(x,x-y)}
{|x-y|^{n-\alpha}} f(y)\hbox{d}y,
\end{equation}
where $\Omega(x,z)\in L^\infty(\rn)\times L^r(S^{n-1}).$ When
$\alpha=0$, the right-hand side integral in (1.1) is interpreted
in the sense of Cauchy principal value. In addition, we assume
that $\Omega(x,z)$ satisfies the cancellation condition,
\begin{equation}
\int_{S^{n-1}}\Omega(x,z')\hbox{d}\sigma(z')=0, \quad \hbox{for
all} \ x\in \rn.\end{equation} Then $T_{\Omega,0}$ is the singular
integral with variable kernel, and we simply write it as
$T_\Omega.$ The $L^p$-boundedness of the singular integral
operator with variable kernel appears in \cite{cz1} (see also
\cite{cz3,cdr}). It turns out that such kind of operators are much
more closely related to the elliptic partial differential
equations of second order with variable coefficients.

\begin{theor}[\cite{cz1}]
Let $r>2(n-1)/n${\rm ,} if $\Omega(x,z)\in L^\infty(\rn)\times
L^r(S^{n-1})$ and satisfies {\rm (1.2),} then there exists a positive
constant $C${\rm ,} independent of $f${\rm ,} such that
\begin{equation*}
\|T_\Omega (f)\|_{L^2(\rn)}\leq C\|f\|_{L^2(\rn)}.
\end{equation*}
\end{theor}

In 1971, Muckenhoupt and Wheeden \cite{mw} established the
following $(L^p, L^q)$-bounded\-ness of $T_{\Omega, \alpha}$ when
$0<\alpha<n$.

\begin{theor}[\cite{mw}]
Let $0<\alpha <n,\ 1<p<n/\alpha$ and $1/q=1/p-\alpha/n$. If there
exists a real number $r>p'$ such that $\Omega \in
L^{\infty}(\rn)\times L^r(S^{n-1})${\rm ,} then there is a positive
constant $C${\rm ,} independent of $f$ and $\Omega${\rm ,} such that
\begin{equation*}
\|T_{\Omega, \alpha}(f)\|_{L^q(\rn)}\leq C\|\Omega\|_
{L^{\infty}(\rn)\times L^r(S^{n-1})}\|f\|_{L^p(\rn)}.
\end{equation*}
\end{theor}

In 2002, Ding {\it et~al} \cite{dcf} studied the boundedness properties
of $\T$ on the Hardy spaces. Recently, Zhang and Ding \cite{zd}
improved some of the results in \cite{dcf} for the case $0<\alpha
<n$, and Chen and Zhang \cite{chz} studied the boundedness
properties of $\T$ on the Herz-type Hardy spaces when $0\leq \alpha
<n$.

The second class of operators considered in this paper are the
Calder\'on--Zygmund singular integrals and the fractional
integrals with homogeneous convolution kernel. Let $\tOmega \in
L^1(S^{n-1})$ be homogeneous of degree 0. We define the integral
operator $T_{\tOmega,\alpha}$ as follows:
\begin{equation}
T_{\tOmega}f(x)=\int_{\rn}\frac{\tOmega(x-y)}{|x-y|^{n-\alpha}}
f(y)\hbox{d}y, \qquad 0\leq \alpha <n.
\end{equation}

When $0<\alpha<n$, $\tT\,$ is the fractional integral, which is
studied by many authors (see \cite{cww,d,dl1,dl3,mw,z} for
instance). We simply write $T_{1,\alpha}$ as $I_{\alpha}$, the
Riesz potential. When $\alpha=0$ the right-hand side integral in
(1.3) is interpreted in the sense of Cauchy principal value. In
addition, we also assume the following cancellation condition:
\begin{equation}
\int_{S^{n-1}}\tOmega(x') \hbox{d}\sigma(x')=0.
\end{equation}
Then $T_{\tOmega,0}$ is the classical Calder\'on--Zygmund singular
integral operator (see \cite{cz2}), and we simply denote it by
$T_{\tOmega}$.

\begin{remark}
{\rm When $\Omega=\tOmega$, Theorem~B is true for $n\geq
n/(n-\alpha)$ (see \cite{m} or \cite{mw} for details).}
\end{remark}

On the other hand, commutators of linear operators take important
roles in harmonic analysis and related topics (see
\cite{clms,crw,j} for example). Let $b$ be a suitable function
such as a BMO or Lipschitz function and $0\leq \alpha<n$. The
commutators generated by $\T$ or $\tT$ and $b$ are defined
formally by
\begin{align}
[b,\T]f(x)&=\int_{\rn}\frac{\Omega(x,x-y)}{|x-y|^{n-\alpha}}
(b(x)-b(y))f(y)\hbox{d}y,\\[.3pc]
[b,\tT]f(x)&=\int_{\rn}\frac{\tOmega(x-y)}{|x-y|^{n-\alpha}}
(b(x)-b(y))f(y)\hbox{d}y.
\end{align}

When $b\in \hbox{BMO}(\rn)$, a study on $[b,T_{\tOmega}]$ has a
long history. In 1976, Coifman {\it et~al} \cite{crw} proved that
$[b,T_{\tOmega}]$ is bounded from $L^p(\rn)\ (1<p<\infty)$ into
itself when $\tOmega \in \hbox{Lip}_1(S^{n-1})$. It is known that
$[b,T_{\tOmega}]$ is neither of weak type (1,1) nor of $(H^1,L^1)$
type (see \cite{pa,pe} for details). When $b\in \hbox{BMO}(\rn)$
and $0<\alpha <n$, the commutator $[b,\tT]$ has been extensively
and profoundly studied by many authors, and we refer the readers
to \cite{ch,dl2,dlz1,dlz2}.

It is well-known that there are other links between the
boundedness properties of commutators and the smoothness of $b$
(see \cite{j,lwz,lwy,lz,pa} for instance). These studies show that
there are much more differences between $b\in \hbox{BMO}(\rn)$ and
$b$ belonging to the Lipschitz space. For $0<\beta \leq 1$, the
Lipschitz space $\L$ is defined by
\begin{equation*}
\L = \left\{g\!\!\,: \|g\|_{{\rm Lip}_{\beta}({\mathbb R}^{n})} =
\dsup_{x,y\in \rn; x\neq y}\frac {|g(x)-g(y)|}{|x-y|^\beta}
<\infty \right\}.
\end{equation*}

Recently, Lu {\it et~al} \cite{lwy} considered the boundedness of
$[b,T_{\tOmega}]$ and $[b, I_{\alpha}]$ on Hardy-type spaces when
$b\in \L$. Motivated by \cite{lwy}, our main aim in this paper is
to study the same problem for commutators $[b,\T]$ and $[b,\tT]$.
We would like to remark that, in \cite{lwy} the smoothness
conditions of $\tOmega \in C^2(S^{n-1})$ or $\tOmega \in
{\rm Lip}_1(S^{n-1})$ are needed, and in this paper we need only a kind
of $L^r$-Dini conditions. So, our results for $[b,\tT]$ improve
and extend the related results in \cite{lwy}, and our computations
are more complex than that in \cite{lwy}. In addition, none of
such kind of results has been seen before for commutators
$[b,\T]$.

Before stating our theorems, we recall the definition of the
$L^r$-Dini condition. For $r\geq 1$, a kernel $\Omega (x,z)$
defined as above satisfies the $L^r$-Dini condition, if
\begin{equation}
\int_0^1 \frac{\omega_r(\delta)}{\delta}\hbox{d}\delta <\infty,
\end{equation}
where $\omega_r(\delta)$ is the integral modulus of order $r$ of
$\Omega$ about $z$, which is defined by
\begin{equation*}
\omega_r(\delta)=\sup_{x\in \rn,|\rho|<\delta}
\left(\int_{S^{n-1}}|\Omega(x,\rho z')-\Omega
(x,z')|^r\hbox{d}\sigma (z') \right)^{1/r},
\end{equation*}
and $\rho$ denotes the rotation in $\rn$ with $|\rho|= \sup_{z'\in
S^{n-1}} |\rho z'-z'|$.

When $\Omega$ does not depend on the first variable, we have
$\Omega=\tOmega$ and we write $\omega_r(\delta)$ as
$\tomega_r(\delta)$.

\setcounter{defin}{0}
\begin{theorem}[\!]
Let $0<\beta \leq 1,\ b\in\L$ and $0\leq
\alpha <n-\beta$. If there is a real number $r>{\rm max}\{n/\beta,
n/(n-\alpha -\beta)\}$ such that $\Omega(x,z)\in
L^{\infty}(\rn)\times L^r(S^{n-1})$ and satisfies the $L^r$-Dini
condition {\rm (1.7),} then there exists a positive constant $C${\rm ,}
independent of $f$ and $b${\rm ,} such that
\begin{equation*}
\|[b,\T]f\|_{L^{n/(n-\alpha-\beta)}(\rn)}\leq C\|b\|_{{\rm Lip}_{\beta}({\mathbb R}^{n})}
\|f\|_{H^1(\rn)}.
\end{equation*}
\end{theorem}

\begin{theorem}[\!]
Let $0<\beta \leq 1,\ b\in \L, \ 0\leq
\alpha <n-\beta$ and $n/(n+\beta)<p<1,\ 1/q = 1/p-(\alpha+\beta)/n$.
If there is a real number $r>{\rm max} \{n/(\beta+n-n/p),
n/(n-\alpha -\beta)\}$ such that $\Omega(x,z) \in L^{\infty}(\rn)
\times L^r(S^{n-1})$ and satisfies
\begin{equation}
\int_0^1\frac{\omega_r(\delta)}{\delta^{1+\beta}}\hbox{\rm d}\delta
<\infty,
\end{equation}
then there exists a positive constant $C${\rm ,} independent of $f$ and
$b${\rm ,} such that
\begin{equation*}
\|[b,\T]f\|_{L^q(\rn)}\leq C\|b\|_{\L}\|f\|_{H^p(\rn)}.
\end{equation*}
\end{theorem}

\setcounter{defin}{1}
\begin{remark}
{\rm As special cases of Theorems~1.1 and 1.2 when
the kernel $\Omega$ does not depend on the first variable, we can
get the corresponding inequalities for $[b,\tT]$.}
\end{remark}

\begin{remark}
{\rm It is worth pointing out that we do not need the cancellation conditions
of $\Omega$ in Theorems~1.1 and 1.2.}
\end{remark}

From Theorems~2.2 and 4.1 in \cite{lwy}, the
$(H^p,L^q)$-boundedness of $[b,\tT]$ fails with $p=n/(n+\beta)$,
even if $\tOmega \in C^2(S^{n-1})$ and satisfies the cancellation
condition (1.4) when $\alpha =0$, so does $[b,I_\alpha]\,(0<\alpha
<n)$. Instead of the $(H^p,L^q)$-boundedness of $[b,\tT]$, we
consider the weak type estimates for $[b,\tT](0\leq \alpha<n)$ when
$p= n/(n+ \beta)$ as in \cite{lwy}. For $0<\alpha<n$, we obtain
the following weak type estimates at the end-point case of
$p=n/(n+\beta)$.

\setcounter{defin}{2}
\begin{theorem}[\!]
Let $0<\beta <1,\ b\in \L$ and $0<
\alpha <n-\beta$. If there is a real number $r\geq n/(n-\alpha)$
such that $\tOmega(x') \in L^r(S^{n-1})$ and satisfies
\begin{equation}
\int_0^1\frac{\tomega_r(\delta)}{\delta^{1+\beta}}{\rm d}\delta<+\infty,
\end{equation}
then $[b,\tT]$ is bounded from $H^{n/(n+\beta)}(\rn)$
into weak $L^{n/(n-\alpha)}(\rn)${\rm ,} that says{\rm ,} there exists a
positive constant $C${\rm ,} independent of $f$ and $b${\rm ,} such that
\begin{align*}
&|\{x\in \rn\!\!: |[b,\tT]f(x)|>\lambda\}|\\[.3pc]
&\quad\,\leq C\|b\|_{\L}
(\lambda^{-1}\|f\|_{H^{n/(n+\beta)}(\rn)})^{n/(n-\alpha)},
\quad \hbox{for all}\ \lambda > 0.
\end{align*}
\end{theorem}

In addition, if we assume $\tOmega$ has integral zero on
$S^{n-1}$, we can establish the following weak type estimate for
$[b,T_{\tOmega}]$ when $p=n/(n+\beta)$.

\begin{theorem}[\!]
Let $0<\beta <1,\ b\in \L$ and $r>1$.
If $\tOmega(x') \in L^r(S^{n-1})$ satisfies {\rm (1.4)} and {\rm (1.9),} then
$[b,T_{\tOmega}]$ is bounded from $H^{n/(n+\beta)}(\rn)$ into weak
$L^1(\rn)${\rm ,} say{\rm ,} there is a positive constant $C${\rm ,} independent of
$f$ and $b${\rm ,} such that
\begin{align*}
&|\{x\in \rn\hbox{\rm :}\ |[b,T_{\tOmega}]f(x)|>\lambda\}|\\[.3pc]
&\quad\,\leq C\|b\|_{\L}
\lambda^{-1}\|f\|_{H^{n/(n+\beta)}(\rn)}, \quad \hbox{for all}\
\lambda>0.
\end{align*}
\end{theorem}

\setcounter{defin}{3}
\begin{remark}
{\rm Theorems~1.3 and 1.4 improve the corresponding results in \cite{lwy}.}
\end{remark}

Our paper is arranged as follows. In $\S 2$, we will formulate
some known results to be used and then prove Theorems~1.1 and 1.2.
In the last section, we prove Theorems~1.3 and 1.4.

\section{Proofs of Theorems~1.1 and 1.2}

\setcounter{equation}{0}

Firstly, we have the following boundedness properties of $[b,\T]$
on $L^p(\rn)$ spaces.

\setcounter{defin}{0}
\begin{proposition}$\left.\right.$\vspace{.5pc}

\noindent Let $0<\beta \leq 1,\ b\in \L,\ 0\leq
\alpha<n-\beta$ and $1<p<n/(\alpha +\beta),\ 1/q=1/p-(\alpha
+\beta)/n$. If $\Omega (x,z)\in L^{\infty}(\rn) \times
L^r(S^{n-1})$ and $r>p'${\rm ,} then there is a positive constant $C${\rm ,}
independent of $f$ and $b${\rm ,} such that\vspace{-.2pc}
\begin{equation*}
\|[b,\T]f\|_{L^q(\rn)}\leq C\|b\|_{\L}\|f\|_{L^p(\rn)}.
\end{equation*}
\end{proposition}

\begin{proof}
By (1.5) and the definition of $\L$, it is easy to
see that\vspace{-.2pc}
\begin{align*}
|[b,\T]f(x)|\leq \|b\|_{\L}T_{|\Omega|, \alpha +\beta}(|f|)(x).
\end{align*}

Applying Theorem~B for $T_{|\Omega|, \alpha +\beta}$, we meet the
desired result.
\hfill $\Box$
\end{proof}

To prove Theorems~1.1 and 1.2, we recall the atomic decomposition
theory of the Hardy spaces. Denote by $[t]$ the greatest integer
which is less than or equal to $t$.

\setcounter{defin}{0}
\begin{defi}$\left.\right.$\vspace{.5pc}

\noindent {\rm Let $0<p\leq 1\leq q<\infty$, $p\neq
q$ and $s\geq [n(1/p-1)]$ be a nonnegative integer. A function
$a(x)\in L^q(\rn)$ is said to be a $(p,q)$-atom centered at $x_0$,
if
\begin{enumerate}
\renewcommand{\labelenumi}{(\alph{enumi})}
\leftskip .1pc
\item supp$(a) \subset B(x_0,d):=\{x\in \rn\!\!: |x-x_0|\leq d\},\ \hbox{for
some}\ d>0$;\vspace{.1pc}

\item $\|a\|_{L^q(\rn)}\leq |B(x_0,d)|^{1/q-1/p};$\vspace{.1pc}

\item $\int a(x)x^{\gamma}{\rm d}x=0,\ \hbox{for} \ 0\leq |\gamma|\leq s$.\vspace{-.5pc}
\end{enumerate}}
\end{defi}

\setcounter{defin}{0}
\begin{lemma}\hskip -.3pc {\rm \cite{lu,st}.}\ Let
$0<p\leq 1$. A distribution $f$ on $\rn$ belongs to $H^p(\rn)$ if
and only if $f$ can be written as $f=\sumj \lambda_j a_j$ in the
sense of distribution{\rm ,} where each $a_j$ is a $(p,q)$-atom{\rm ,}
$\lambda_j \in {\Bbb C}$ and $\sumj |\lambda_j|^p <\infty$.
Furthermore{\rm ,}\vspace{-.2pc}
\begin{equation*}
\|f\|_{H^p(\rn)} \sim \inf \left\{\sumj
|\lambda_j|^p\right\}^{1/p},
\end{equation*}
where the infimum is taken over all the above atomic
decompositions of $f$.
\end{lemma}

Now, we also need the following estimate for the kernel
$\Omega(x,z)$, which is stated in \cite{dcf} and can be proven by
using the same methods as that of Lemma~5 in \cite{kw}.

\begin{lemma}\hskip -.3pc{\rm \cite{dcf}.}\
Let $0\leq \mu <n${\rm ,}
$\Omega(x,z) \in L^\infty(\rn)\times L^r(S^{n-1})\ (r\geq 1)$ and
satisfies the $L^r$-Dini condition {\rm (1.7)}. If there is a constant
$0< a_0<1/2$ such that $|y|<a_0R${\rm ,} then
\begin{align*}
&\left(\dint_{R<|x|<2R} \left|\dfrac {\Omega(x, x-y)}
{|x-y|^{n-\mu}}-\dfrac{\Omega(x, x)}{|x|^{n-\mu}}
\right|^r{\rm d}x\right)^{1/r}\\[.3pc]
&\quad \leq CR^{n/r-(n-\mu)}\left\{\dfrac{|y|}{R}+\dint_{|y|/2R <\delta
<|y|/R}\frac{\omega_r(\delta )}{\delta} {\rm d}\delta\right\},
\end{align*}
where $C$ is a positive constant independent of $R$ and $y$.
\end{lemma}

\setcounter{section}{1}
\setcounter{defin}{0}
\begin{pot}
{\rm Since $r>n/(n-\alpha -\beta)$,
$1\leq r'<n/(\alpha +\beta)$. Pick a real number $\ell_1$ such that
$r'<\ell_1 <n/(\alpha +\beta)$ and set $1/\ell_2=1/\ell_1
-(\alpha+\beta)/n$. From Proposition~2.1 we see that $[b,\T]$ is
bounded from $L^{\ell_1}(\rn)$ into $L^{\ell_2}(\rn)$.

Write $q=n/(n-\alpha-\beta)$ for simplicity. By the atomic
decomposition theory of Hardy space, to prove Theorem~1.1, it
suffices to verify that for any $(1,\ell_1)$-atom $a(x)$, there is
a constant $C>0$, independent of $a$ and $b$, such that
$\|[b,\T]a\|_{L^q(\rn)}\leq C\|b\|_{\L}$. Without loss of
generality, we assume supp$(a)\subset B=B(0,d)$. Denote by
$tB=B(0,td)$ for $t>0$, then
\setcounter{section}{2}
\setcounter{equation}{0}
\begin{align}
\|[b,\T]a\|_{L^q(\rn)}&\leq \left(\dint_{2B}|[b,\T]a(x)|^q\hbox{d}x
\right)^{1/q}\nonumber\\[.3pc]
&\quad\, + \left(\dint_{(2B)^C}|[b,\T]a(x)|^q\hbox{d}x\right)^{1/q} := \hbox{I} + \hbox{II}.
\end{align}

Noting that $\ell_2 >n/(n-\alpha-\beta)=q$, by the H\"older's
inequality, the $(L^{\ell_1},L^{\ell_2})$-boundedness of $[b,\T]$
and the size condition of $a$, we have
\begin{align}
\hbox{I} &\leq |2B|^{1/q-1/\ell_2}\|[b,\T]a\|_{L^{\ell_2}(\rn)}\nonumber\\[.2pc]
&\leq C\|b\|_{\L}|B|^{1/q-1/\ell_2}\|a\|_{L^{\ell_1}(\rn)}\nonumber\\[.2pc]
&\leq C\|b\|_{\L}.
\end{align}
By the cancellation condition of $a$ and the Minkowski's
inequality, we have
\begin{align}
\hbox{II} &\leq \left(\dint_{(2B)^C}\left|(b(x)-b(0))\dint_B
\left[\dfrac {\Omega(x,x-y)}{|x-y|^{n-\alpha}}-\dfrac{\Omega(x,x)}
{|x|^{n-\alpha}}\right]a(y)\hbox{d}y\right|^q \hbox{d}x \right)^{1/q}\nonumber\\[.2pc]
&\quad\, + \left(\dint_{(2B)^C}\left|\dint_B\dfrac{\Omega(x,x-y)}
{|x-y|^{n-\alpha}}\left(b(y)-b(0)\right)a(y)\hbox{d}y\right|^q\hbox{d}x\right)^{1/q}\nonumber\\[.2pc]
&:= \hbox{II}_1 + \hbox{II}_2.
\end{align}

Denote by $\Delta(x,y)=\frac{\Omega(x,x-y)}{|x-y|^{n-\alpha}}
-\frac{\Omega(x,x)}{|x|^{n-\alpha}}$ for simplicity. Note that
$r>n/(n-\alpha-\beta)=q$, by the H\"older's inequality and Lemma~2.2,
for any $y\in B$,
\begin{align}
&\left(\dint_{2^jd<|x|<2^{j+1}d}\left|\Delta(x,y)\right|^q\hbox{d}x\right)^{1/q}\nonumber\\
&\quad\,\leq(2^jd)^{n(1/q-1/r)}\left(\dint_{2^jd<|x|<2^{j+1}d}
|\Delta(x,y)|^r\hbox{d}x\right)^{1/r}\nonumber\\[.2pc]
&\quad\,\leq C(2^jd)^{n/q-n+\alpha}\left\{2^{-j}
+\dint_{|y|/2^{j+1}d}^{|y|/2^jd}\dfrac{\omega_r(\delta)}
{\delta}\hbox{d}\delta\right\}.
\end{align}

Noting that $\sup_{x\in2^{j+1}B}|b(x)-b(0)|\leq
\|b\|_{\L} (2^{j+1}d)^{\beta}$ and $q=n/(n-\alpha-\beta)$, by
(2.4), (1.7) and the size condition of $a$, we have
\begin{align*}
\hskip -4pc \hbox{II}_1&\leq C\dint_B|a(y)|\dsum_{j=1}^{\infty}\left(\dint_{2^jd<|x|<
2^{j+1}d}|\Delta(x,y)|^q|b(x)-b(0)|^q\hbox{d}x\right)^{1/q}\hbox{d}y\\[.2pc]
&\leq C\|b\|_{\L}\dint_B|a(y)| \dsum_{j=1}^{\infty} (2^jd)^{\beta}
\left(\dint_{2^jd<|x|<2^{j+1}d}|\Delta(x,y)|^q\hbox{d}x\right)^{1/q}\hbox{d}y
\end{align*}
\begin{align}
&\leq C\|b\|_{\L}\dint_B|a(y)|\dsum_{j=1}^{\infty}
\left\{2^{-j}+\dint_{|y|/2^{j+1}d}^{|y|/2^jd}
\dfrac{\omega_r(\delta)}{\delta}\hbox{d}\delta\right\}\hbox{d}y\nonumber\\[.2pc]
&\leq C\|b\|_{\L}\dint_B|a(y)|\hbox{d}y \leq C\|b\|_{\L}.
\end{align}

To consider $\hbox{II}_2$, for any fixed $x\in (2B)^C$, we first estimate
$\int_B|\Omega(x,x-y)||a(y)|\hbox{d}y$. Noting that $r'<\ell_1$ and for
any $x\in (2B)^C$ and $y\in B$ there holds $|x|/2<|x-y|<2|x|$,
then by the H\"older's inequality, we have for any fixed $x\in
(2B)^C$,
\begin{align}
\dint_B&|\Omega(x,x-y)||a(y)|\hbox{d}y\leq \left(\dint_B|\Omega(x,x-y)|^r
\hbox{d}y\right)^{1/r}\left(\dint_B|a(y)|^{r'}\hbox{d}y\right)^{1/r'}\nonumber\\
&\leq C\|a\|_{L^{\ell_1}(\rn)}|B|^{1/r'-1/\ell_1} \left(
\dint_{|x|/2<|x-y|<2|x|}|\Omega(x,x-y)|^r\hbox{d}y\right)^{1/r}\nonumber\\
&\leq C\|\Omega\|_{L^{\infty}(\rn)\times L^r(S^{n-1})}
|B|^{1/r'-1/\ell_1}\|a\|_{L^{\ell_1}(\rn)}|x|^{n/r}.
\end{align}

Since $r>n/\beta$ and $q=n/(n-\alpha-\beta)$, then
$(n/r-n+\alpha)q=q(n/r-\beta)-n<-n$. By the size condition of $a$,
we have
\begin{align}
\hbox{II}_2&\leq C\|b\|_{\L}|B|^{\beta/n}\!\left\{\!\dint_{(2B)^C}\!\left(\!
|x|^{-n+\alpha}\!\dint_B|\Omega(x,x\!-\!y)||a(y)|\hbox{d}y\!\right)^q \!\hbox{d}x\!\right\}^{\!1/q}\nonumber\\
&\leq C\|b\|_{\L}|B|^{\beta/n+1/r'-1/\ell_1} \|a\|_{L^{\ell_1}
(\rn)}\left(\dint_{(2B)^C}|x|^{(n/r-n+\alpha)q}\hbox{d}x\right)^{1/q}\nonumber\\
&\leq C\|b\|_{\L}|B|^{\beta/n-1/r} \left(\dint_{2d}^{\infty}
{\rho}^{q(n/r-\beta)-1}\hbox{d}\rho\right)^{1/q}\nonumber\\
&\leq C\|b\|_{\L}.
\end{align}

The above estimates for $\hbox{II}_1$ and $\hbox{II}_2$, together with
(2.3) show that $\hbox{II}\leq C\|b\|_{\L}$. This completes the proof of
Theorem~1.1.}\hfill $\Box$
\end{pot}

\setcounter{section}{1}
\begin{pot}
{\rm Similar to the proof of Theorem~1.1, we assume $\ell_1$ and
$\ell_2$ to be the same as in the proof of Theorem~1.1. For any
$(p,\ell_1)$-atom $a$ with supp$(a)\subset B(0,d)$, it suffices to
verify that there is a constant $C>0$, independent of $a$ and $b$, such
that $\|[b,\T]a\|_{L^q(\rn)}\leq C\|b\|_{\L}$. Write
\begin{equation*}
\|[b,\T]a\|_{L^q(\rn)}\leq \hbox{I} + \hbox{II}_1 + \hbox{II}_2,
\end{equation*}
where $\hbox{I},\ \hbox{II}_1$ and $\hbox{II}_2$ are the same as in (2.1) and (2.3),
except for $n/(n+\beta)<p<1$ and $1/q=1/p-(\alpha+ \beta)/n$.

Since $0<p<1<r'<\ell_1$ then $q<\ell_2$. Similar to (2.2) we have
$\hbox{I}\leq C\|b\|_{\L}$. To finish the proof of Theorem~1.2, we need
only to modify the estimates for $\hbox{II}_1$ and $\hbox{II}_2$ in the proof of
Theorem~1.1.

Noting that $r>n/(n-\alpha -\beta),\ n/(n+\beta)<p<1$ and $0<\beta
\leq 1$, it is easy to see that $r>q$ and $n/p-n-1\leq
n/p-n-\beta<0$. Since (2.4) is always true for $r>q$, then by
(1.8) and the size condition of $a$, similar to (2.5), we have
\pagebreak

$\left.\right.$\vspace{-1.5pc}
\begin{align*}
\hskip -4pc \hbox{II}_1&\leq \left(\dint_{(2B)^C}\left|(b(x)-b(0))\dint_B\left[\dfrac
{\Omega(x,x-y)}{|x-y|^{n-\alpha}} -\dfrac {\Omega(x,x)}
{|x|^{n-\alpha}}\right]a(y)\hbox{d}y \right|^q \hbox{d}x\right)^{1/q}\\[.3pc]
\hskip -4pc &\leq C\|b\|_{\L}\dint_B|a(y)|\dsum_{j=1}^{\infty}
(2^jd)^{n/p-n}\left\{2^{-j}+\dint_{|y|/2^{j+1}d}^{|y|/2^jd}
\dfrac{\omega_r(\delta)}{\delta}\right\}\hbox{d}y\\[.3pc]
\hskip -4pc &\leq C\|b\|_{\L}d^{n/p-n}\dint_B|a(y)|\dsum_{j=1}^{\infty}
2^{j(n/p-n)}\left\{2^{-j}+2^{-j\beta}\dint_0^1
\dfrac{\omega_r(\delta)}{\delta^{1+\beta}}\hbox{d}\delta\right\}\hbox{d}y\\[.3pc]
\hskip -4pc &\leq C\|b\|_{\L}d^{n/p-n}\dint_B|a(y)|\hbox{d} y\left(1+\dint_0^1
\dfrac{\omega_r(\delta)}{\delta^{1+\beta}}\hbox{d}\delta\right)\\[.3pc]
\hskip -4pc &\leq C\|b\|_{\L}.
\end{align*}

From $r>n/(\beta +n -n/p)$ and $1/p=1/q+(\alpha +\beta)/n$, we
have $n+q(n/r-n+\alpha)=q(n/r-\beta-n+n/p)<0$. By (2.6), similar
to (2.7), we get
\begin{align*}
\hskip -4pc \hbox{II}_2&\leq
C\|b\|_{\L}|B|^{\beta/n}\!\left\{\dint_{(2B)^C}\!|x|^{-q(n-\alpha)}\!\left(
\dint_B|\Omega(x,x-y)||a(y)|\hbox{d}y\right)^q\!\hbox{d}x\!\right\}^{1/q}\\[.3pc]
\hskip -4pc &\leq C\|b\|_{\L}|B|^{\beta/n+1/r'-1/\ell_1}
\|a\|_{L^{\ell_1}(\rn)}\left(\dint_{(2B)^C}
|x|^{q(n/r-n+\alpha)}\hbox{d}x\right)^{1/q}\\[.3pc]
\hskip -4pc &\leq C\|b\|_{\L}|B|^{(\alpha +\beta)/n-(1/p-1/q)}\leq C\|b\|_{\L}.
\end{align*}
Summing up the discussion above, we finish the proof of
Theorem~1.2.}\hfill $\Box$\vspace{-.4pc}
\end{pot}

\setcounter{section}{2}
\section{Proofs of Theorems~1.3 and 1.4}
\setcounter{equation}{0}
\setcounter{defin}{0}

In this section, we prove Theorems~1.3 and 1.4. To do this, we
need the following known estimates for  $\tT \,(0<\alpha <n)$ and
$T_{\tOmega}$.

\begin{lemma}\hskip -.3pc{\rm \cite{cww,d}.}\
Let $0<\alpha<n, \
r\geq n/(n-\alpha)$ and $\tOmega \in L^r(S^{n-1})${\rm ,} then $\tT$ is
of weak type $(1,n/(n-\alpha))${\rm ,} that is{\rm ,} there is a positive
constant $C${\rm ,} such that
\begin{align*}
|\{x\in \rn\hbox{\rm :}\ |\tT f(x)|>\lambda\}|\leq C(\lambda^{-1}
\|f\|_{L^1(\rn)})^{n/(n-\alpha)},\quad \hbox{for all}\ \lambda>0.
\end{align*}
\end{lemma}

\begin{lemma}\hskip -.3pc{\rm \cite{cz2,s}.}\
Suppose that $\tOmega \in L\log L (S^{n-1})$ and satisfies {\rm (1.4)}. Then
$T_{\tOmega}$ extends to an operator of type $(p,p)$ for
$1<p<\infty${\rm ,} and of weak type $(1,1)$.
\end{lemma}

\setcounter{section}{1}
\setcounter{defin}{2}
\begin{pot}
{\rm Denote by $q_0=n/(n-\alpha)$ for
simplicity. For $f\in H^{n/(n+\beta)}(\rn)$, by the atomic
decomposition theory of Hardy space, $f=\sum_{j=-\infty}^{\infty}
\lambda_j a_j$ in the sense of distribution, where each $a_j$ is a
$(n/(n+\beta),\ell_1)$-atom and $\lambda_j\in {\Bbb C}$. Suppose
that supp$(a_j)\subset B_j=B(x_j,r_j)$, then
\begin{align*}
[b,\tT]f(x)&=\dsumj \lambda_j (b(x)-b(x_j))\tT
a_j(x)\chi_{2B_j}(x)
\end{align*}
\begin{align*}
&\quad\, +\dsumj \lambda_j (b(x)-b(x_j))\tT
a_j(x)\chi_{(2B_j)^C}(x)\\
&\quad\, - \tT \left(\dsumj \lambda_j (b-b(x_j))a_j\right)(x)\\[.3pc]
&:= \hbox{I}_1(x) + \hbox{I}_2(x) + \hbox{I}_3(x).
\end{align*}
\setcounter{section}{3}
From Remark~1.1, we can choose $\ell_1$ and $\ell_2$ with
$1<\ell_1<n/\alpha$ and $1/\ell_2=1/\ell_1 -\alpha/n$ such that
$\tT$ is bounded from $L^{\ell_1}(\rn)$ into $L^{\ell_2}(\rn)$.
Noting that $b\in \L$ and $\ell_2 >q_0$, by the size condition of
$a_j$, we get
\begin{align*}
&\|(b-b(x_j))\tT a_j \chi_{2B_j}\|_{L^{q_0}(\rn)}\\
&\qquad \leq C\|b\|_{\L}|B_j|^{\beta/n}\|\tT a_j\chi_{2B_j}\|_{L^{q_0}(\rn)}\\[.3pc]
&\qquad \leq C\|b\|_{\L}|B_j|^{\beta/n+1/q_0 -1/\ell_2} \|\tT
a_j\|_{L^{\ell_2}(\rn)}\\[.3pc]
&\qquad \leq C\|b\|_{\L}|B_j|^{\beta/n+1/q_0-1/\ell_2}
\|a_j\|_{L^{\ell_1}(\rn)}\\[.3pc]
&\qquad \leq C\|b\|_{\L}.
\end{align*}
Consequently,
\begin{align}
|\{x\in \rn\hbox{\rm :}\ &|{\rm I}_1(x)|>\lambda/3\}|^{1/q_0}\leq
3\lambda^{-1}\|{\rm I}_1\|_{L^{q_0}(\rn)}\nonumber\\[.3pc]
&\le 3\lambda^{-1}\dsumj |\lambda_j|\|(b-b(x_j))\tT
a_j\chi_{2B_j}\|_{L^{q_0}(\rn)}\nonumber\\[.3pc]
&\leq C\|b\|_{\L}\lambda^{-1}\dsumj |\lambda_j|.
\end{align}

Noting that $b\in \L$, by the H\"older's inequality and the size
condition of $a_j$,
\begin{align*}
\|(b-b(x_j))a_j\|_{L^1(\rn)}&\leq
\|b\|_{\L}|B_j|^{\beta/n}\|a_j\|_{L^1(\rn)}\\[.3pc]
&\leq C\|b\|_{\L}|B_j|^{\beta/n +1-1/\ell_1}
\|a_j\|_{L^{\ell_1}(\rn)}\\[.3pc]
&\leq C\|b\|_{\L}.
\end{align*}

From Lemma~3.1, $\tT$ is of weak type $(1, q_0)$, then we have
\begin{align}
|\{x\in \rn\hbox{:}\ |\hbox{I}_3(x)|>\lambda/3\}|^{1/q_0}&\leq 3 \lambda^{-1}\dsumj
|\lambda_j|\|(b-b(x_j))a_j\|_{L^1(\rn)}\nonumber\\[.3pc]
&\leq C\|b\|_{\L}\lambda^{-1} \dsumj |\lambda_j|.
\end{align}

Now, we are in position to give the same estimates for $\hbox{I}_2(x)$ as
that of $\hbox{I}_1(x)$ and $\hbox{I}_3(x)$. Write $\Delta
(x,y,x_j)=\frac{\tOmega(x-y)}{|x-y|^{n-\alpha}} -\frac
{\tOmega(x-x_j)}{|x-x_j|^{n-\alpha}}$ for simplicity. If
$r>n/(n-\alpha)=q_0$, then by using the H\"older's inequality and
Lemma~2.2 for the special case that the kernel does not depend on
the first variable, we have, for any $y\in B_j$ and $k\geq 1$,
\begin{align*}
&\left(\dint_{2^kr_j<|x|<2^{k+1}r_j} |\Delta
(x,y,x_j)|^{q_0}\hbox{d}x\right)^{1/q_0}\\[.2pc]
&\qquad \leq C(2^kr_j)^{n-\alpha -n/r}\left(\dint_{2^kr_j<|x|<2^{k+1}r_j}
|\Delta(x,y,x_j)|^r\hbox{d}x\right)^{1/r}\\[.2pc]
&\qquad \leq C\left\{2^{-k}+ \dint_{|y|/2^{k+1}r_j}^{|y|/2^kr_j}
\dfrac{\tomega_r(\delta)}{\delta}\hbox{d}\delta\right\}.
\end{align*}
Obviously, there holds the same estimate when
$r=n/(n-\alpha)=q_0$.

By the cancellation condition of $a_j$, the Minkowski's
inequality, the above estimate and (1.9), we have
\begin{align*}
\hskip -4pc &\|(b-b(x_j))\tT (a_j)\chi_{(2B_j)^C}\|_{L^{q_0}(\rn)}\\[.2pc]
\hskip -4pc &\quad\,\leq \dint_{B_j}|a_j(y)|\dsum_{k=1}^{\infty}\left(\dint_{2^kr_j<
|x|<2^{k+1}r_j}|\Delta (x,y,x_j)(b(x)-b(x_j))|^{q_0}
\hbox{d}x\right)^{1/q_0}\hbox{d}y\\[.2pc]
\hskip -4pc &\quad\,\leq C\|b\|_{\L}\!\dint_{B_j}\!|a_j(y)|\dsum_{k=1}^{\infty}
(2^kr_j)^{\beta}\!\left(\dint_{2^kr_j<|x|<2^{k+1}r_j}\!|\Delta
(x,y,x_j)|^{q_0}\hbox{d}x\!\right)^{1/q_0}\!\hbox{d}y\\[.2pc]
\hskip -4pc &\quad\,\leq C\|b\|_{\L}r_j^{\beta}\dint_{B_j}|a_j(y)|
\dsum_{k=1}^{\infty}2^{k\beta}\left\{2^{-k}+
\dint_{|y|/2^{k+1}r_j}^{|y|/2^kr_j} \dfrac{\tomega_r(\delta)}
{\delta}\hbox{d}\delta\right\}\hbox{d}y\\[.2pc]
\hskip -4pc &\quad\,\leq C\|b\|_{\L}r_j^{\beta}\dint_{B_j}|a_j(y)|
\dsum_{k=1}^{\infty}\left\{2^{k(\beta-1)}+
\dint_{|y|/2^{k+1}r_j}^{|y|/2^kr_j}\dfrac{\tomega_r(\delta)}
{{\delta}^{1+\beta}}\hbox{d}\delta\right\}\hbox{d}y\\[.2pc]
\hskip -4pc &\quad\,\leq C\|b\|_{\L}|B_j|^{\beta/n+1-1/\ell_1}
\|a_j\|_{L^{\ell_1}(\rn)} \leq C \|b\|_{\L}.
\end{align*}
And then\vspace{-.2pc}
\begin{align}
|\{x\in \rn\hbox{\rm :}\ &|\hbox{I}_2(x)|>\lambda/3\}|^{1/q_0}\leq
C\lambda^{-1}\|\hbox{I}_2\|_{L^{q_0}(\rn)}\nonumber\\[.2pc]
&\leq C\lambda^{-1}\dsumj|\lambda_j|\|(b-b(x_j))
\tT(a_j)\chi_{2B_j}\|_{L^{q_0}(\rn)}\nonumber\\
&\leq C\|b\|_{\L}\lambda^{-1}\dsumj |\lambda_j|.
\end{align}
From (3.1)--(3.3), and noting that $n/(n+\beta)<1$, we
have\vspace{-.2pc}
\begin{align*}
\hskip -4pc |\{x\in \rn\hbox{\rm :}\ |[b,\tT]|>\lambda\}|^{1/q_0}&\leq \dsum_{i=1}^3| \{x\in
\rn\!\!:|I_i(x)|>\lambda/3\}|^{1/q_0}\\
\hskip -4pc &\leq C\|b\|_{\L}\lambda^{-1}\left(\dsumj |\lambda_j|^{n/(n+\beta)}
\right)^{(n+\beta)/n}.
\end{align*}

\noindent This completes the proof of Theorem~1.3.}\hfill$\Box$
\end{pot}
\pagebreak

\setcounter{section}{1}
\begin{pot}
{\rm In the proof of Theorem~1.3, set
$\ell_1=\ell_2>1$, applying Lemma~2.2 for $\mu =0$ and Lemma~3.2,
we can obtain the desired result, we omit the details.} \hfill $\Box$
\end{pot}

\section*{Acknowledgement}

The authors would like to thank the
referees for their comments and suggestions. This project is
supported by the Research Funds of Zhejiang Sci-Tech University
(No. 0313055-Y) and NSFZJ (No. Y604563).

\end{document}